# On the Convex Feasibility Problem

**Laura Măruşter[1] and Ştefan Măruşter[2]**
[1] University of Gröningen, Gröningen, Nederland
l.maruster@rug.nl
[2] West University of Timişoara, Romania
maruster@info.uvt.ro

**ABSTRACT**. The convergence of the projection algorithm for solving the convex feasibility problem for a family of closed convex sets, is in connection with the regularity properties of the family. In the paper [18] are pointed out four cases of such a family depending of the two characteristics: the emptiness and boudedness of the intersection of the family. The case four (the interior of the intersection is empty and the intersection itself is bounded) is unsolved. In this paper we give a (partial) answer for the case four: in the case of two closed convex sets in $\mathbb{R}^3$ the regularity property holds.

## Introduction

The convex feasibility problem can be formulated in a very simple way: Let $A_i$, $i=1,...,m$, be a family of closed convex sets with nonempty intersection, $\cap A_i \neq \emptyset$, in a real Hilbert space; the convex feasibility problem is:

*Find a point in $\cap A_i$.*

Historically, this problem arisen in connection with guessing a starting point in SIMPLEX algorithm for solving the mathematical programming problem. In this particular case, the sets of family are given by a finite number of halfspaces defined by linear inequalities; a point which satisfies all these inequalities is a "feasible point" and it can be used for starting the iteration

225

process. A first result for solving the convex feasibility problem in this case was given by Agmon, Motzkin and Schoenberg in 1954, [1, 20].

In the last decades a number of papers were written on the subject and results concerning some theoretical aspects and particularly algorithms for solving this problem were obtained. The interest for the convex feasibility problem is motivated by some important applications in specific areas. In the paper [2] are pointed out the main applications both in other mathematical algorithms and directly in some practical problems; we list here a part of them.

*Best approximation* with applications in linear prediction theory, partial differential equations (Dirichlet problem), complex analysis (Bergman kernels, conformal mappings);

*Subgradient algorithms* with application in solution of convex inequalities, minimization of convex nonsmooth functions;

*Image reconstruction, discrete and continuous models* with applications in radiation therapy treatment planing, electron microscopy, computerized tomography, signal processing.

Usually, the convex feasibility problems are solved by *projection algorithms*. The geometric idea of the projection method is to project the current iteration onto certain set form the intersecting family and to take the next iteration on the straight line connecting the current iteration and this projection. A weight factor gives the exact position of the next iteration. Different strategies concerning the selection of the set onto which the current iteration will be projected, will give particular projection-type algorithms.

If $P_{M_i}(x)$ denote the projection of $x$ onto $M_i$, then the classical projection method is

$$x_{k+1} = (1-t_k)x_k + t_k P_{M_{\alpha(k)}}(x_k) \qquad (1)$$

where $t_k$ is the weight factor, $0 < t_k < 2$, and the function $\alpha : \mathbb{N} \to \{1,\ldots,N\}$ defines the strategy. The usual strategy is the cyclic covering of the sets of the family, that is $\alpha(k) = mod_N(k) + 1$

The projection algorithm was used in [1, 20] for solving a system of linear inequalities (the authors referred to their method as *relaxation algorithm*). Generalizations for convex sets in real *n*-dimensional spaces were given in [8,11]. Bergman [4] considered the classical projection method for the case of *m* intersecting closed convex sets $M_i$ in a real Hilbert space. He showed that, given an arbitrary starting point $x_0$, the sequence generated by the projection algorithm converges weakly to a point in $M = \bigcap_{i=1}^{m} M_i$. A complete and exhaustive study on algorithms for solving convex feasibility problem, including comments about their applications and an excellent bibliography, was given by Bauschke and Borwein [2].

The paper is organized as follows. In section 2 we briefly state the convergence properties of the Mann iteration process. The projection algorithm is a particular case of this iteration, so that its convergence results from general theorems of the Mann iteration. Section 3 is devoted to regularity properties of a family of closed convex sets, properties which play a significant roll of the behavior of the sequence generated by this algorithm.

## 1 The projection algorithms and normal Mann iteration

The projection method is a particular case of the Mann iteration process:

$$x_{k+1} = (1-t_k)x_k + t_k T(x_k) \qquad (2)$$

where $t_k$ is a sequence of real numbers satisfying some properties, usually called the control sequence. The convergence properties of the projection algorithm for convex feasibility problem are obtained from the general convergence properties of the Mann iteration.

*Remark 1.* Note that (2) is a particular case of the general Mann iteration $x_{k+1} = T(\bar{x}_k)$, where $\bar{x}_k = \sum_{j=0}^{k} \alpha_{kj} x_j$ and $A = \{\alpha_{kj}\}$ is a triangular averaging matrix. If this matrix satisfies the s*egmenting condition*, that is, $\alpha_{n+1,j} = (1-\alpha_{n+1,n+1})\alpha_{nj}$, then the general Mann iteration becomes just (2) with a specific relaxation strategy, $t_k = \alpha_{k+1,k+1}, \forall k \in N$. If $0 < t_k < 1$, then $x_{k+1}$ is a convex combination of $x_k$ and $T(x_k)$. This restriction concerning $\{t_k\}$ is not always satisfied; the typical case is, for example, the projection algorithm for convex feasibility problem, algorithm which has the form (2) with $0 < t_k < 2$. In particular, if $t_k = \frac{1}{2}$, (2) becomes $x_{k+1}=(x_k+T(x_k))/2$, which is the well known Krasnoselski method. The term *Krasnoselski/Mann* or *the relaxed* iteration is sometimes used for (1.2) as well.

The convergence properties of the normal Mann iteration are in connection with some variational properties of the mapping *T*. Let $\mathcal{H}$ be a real Hilbert space endowed with scalar product $\langle .,. \rangle$ and norm $\|.\|$ and let *C* be a closed convex set in $\mathcal{H}$. If $T: C \to C$ is a nonlinear mapping, *Fix(T)*

227

will denote the set of fixed point of *T* in *C*, set which we will suppose throughout this paper to be nonempty.

**Definition 1**. *The mapping T is said to be qusi-nonexpansive if*
$$\|T(x)-x^*\| \le \|x-x^*\|, \quad \forall x \in C, \quad x^* \in Fix(T)$$

**Definition 2.** *The mapping T is said to be demicontractive (or **k**-demi-contractive) if there exists* **k** < *1 such that*
$$\|T(x)-x^*\|^2 \le \|x-x^*\|^2 + \mathbf{k}\|T(x)-x\|^2, \quad \forall x \in C, \quad x^* \in Fix(t) \qquad (3)$$

Obviously, the class of demicontractive mappings properly includes the class of quasi-nonexpansive mappings for $0 \le \mathbf{k} \le 1$

*Remark 2.* For negative values of **k** the class of demicontractive mappings is diminished in a great extent; in [2] such a class (with negative value of **k**) was considered under the name of *strongly attracting map*. In particular, the mapping *T* which satisfied (3) with **k** = -1 is called *pseudo-contractive* in [26]. Note also that a mapping *T* satisfying (3) with **k** = 1 is usually called *hemicontractive* and it was considered by some authors in connection with strong convergence of the implicit Mann-type iteration (see, for example, [23]).

The notion of quasi-nonexpansivity was introduced by Tricomi in 1916 [25] for a real function *f* defined on a finite or infinite interval (*a,b*) with the values in the same interval. He proved that the sequence $\{x_k\}$ generated by the simple iteration $x_{k+1}=f(x_k)$, $x_0$-given in (*a,b*), converges to a fixed point of *f* provided that *f* is continuous and strictly quasi-nonexpansive on (*a,b*). Stepleman in a paper published in 1975 [24] studied necessary and sufficient conditions for the convergence of this sequence in real case. His main result states that convergence is assured if and only if the second iterate of *f* is strictly quasi-nonexpansive. The importance of the concept of quasi-nonexpansivity for the computation of fixed points in more general cases had been emphasized by many authors (see, for example, the survey papers of Petryshyn and Williamson [22] and of Diaz and Metcalf [7]) and this class of mappings is still being studied extensively (see, for instance, the recent monographs of Chidume [5] and Berinde [3] and the references therein).

The condition of demicontractivity or the more restrictive condition of quasi-nonexpansivity is not sufficient for the convergence of Mann iteration, even in finite dimensional spaces; some additional smoothness properties of the mapping *T*, like continuity or demiclosedness are required.

**Definition 3.** A mapping $T$ is said to be demiclosed at $y$, if for any sequence $\{x_k\}$ which converges weakly to $x$, and if the sequence $\{T(x_k)\}$ converges strongly to $y$, then $T(x)=y$.

In the sequel, as often as not, only the demiclosednes at 0 is used, which is a particular case when $y=0$.

The class of mappings satisfying the condition (3) and the name of *demicontractive* were introduced by Hicks and Kubicek in 1977 [10]. They studied the convergence properties of a sequence $\{x_k\}$ generated by the Mann-type iteration to a fixed point of $T$ in real Hilbert spaces. They proved that if $T$ is demicontractive and if $I-T$ is demiclosed at zero, then the sequence $\{x_k\}$ generated by the Mann iteration (2) converges weakly to a fixed point of $T$. The control sequence is assumed to satisfy the condition $t_k \to t$, $0 < t < 1 - \mathbf{k}$.

In [15] a class of mappings which satisfies so called *condition (A)*:

$$\langle x - T(x), x - x^* \rangle \geq \lambda \|x - T(x)\|^2, \quad \forall x \in C, \quad x^* \in Fix(T) \tag{4}$$

where $\lambda$ is a positive number was considered. It is routine to see that the conditions (3) and (4) are equivalent $\lambda = (1-\mathbf{k})/2$ (indeed, it can be simply checked that $\|x - x^*\|^2 + \mathbf{k}\|x - T(x)\|^2 - \|T(x) - x^*\|^2 = 2\langle x - x^*, x - T(x) \rangle - (1-\mathbf{k})\|x - T(x)\|^2$). Thus the class of demicontractive mappings coincides with the class of mappings satisfying the condition (A). In [15], the same result concerning the weak convergence of the normal Mann iteration was obtained, more exactly, if $T$ satisfies the condition (A) and $I-T$ is demiclosed at zero, then the sequence $\{x_k\}$ converges weakly to a fixed point. The control sequence satisfies a similar condition (to a certain extent, weaker) $0 < a \leq t_k \leq b < 2\lambda$ (or $0 < a \leq t_k \leq b < 1 - \mathbf{k}$). Note that the equiva-lence between the conditions (3) and (4) was observed by some authors [12-14, 6, 21, 19]. Earlier, in 1973, in the paper [16], a similar result in a finite dimensional spaces was presented.

Petryshyn and Williamson [22] pointed out the significant role of the behavior of a sequence with respect to the set of fixed points. For strong convergence of the Mann iteration, it seems that such property plays a prominent part, so that in [17] we have suggested to capture it in a definition.

**Definition 4.** *Let $\{x_k\}$ be a sequence in $\mathcal{H}$ and let $M \subset \mathcal{H}$ be a closed subset. We say that $\{x_k\}$ is* regular with respect *to $M$ if $d(x_k, M) \to 0$ as $k \to \infty$*

**Theorem 1.** *(Petryshyn and Williamson) Suppose that $T : D \subset \mathcal{H} \to \mathcal{H}$ is a quasi-nonexpansive mapping and that Fix(T) is nonempty and closed. Let $x_0 \in D$ such that $x_k = T(x_0)^k \in D, k = 1,2,\ldots$ Then the sequence $\{x_k\}$ converges (strongly) to a fixed point of $T$ if and only if $\{x_k\}$ is regular with respect to Fix(T).*

Here, as usual, $T^k$ denotes the $k^{th}$ iterate of $T$.

*Remark 3.* Theorem (1) is a slight generalization of the first result of [22] and its proof is, practically, identical. Essentially, Theorem (1) replaced the condition of continuity of $T$, from the original result, with the condition of closedness of $Fix(T)$. It is easy to see that the latter condition is weaker, and, as will be seen, is essential for our development.

Consider now the following strategy in the projection algorithm. Let $i_x$ be the least index such that

$$\|x - P(x, i_x)\| = \max_i \|x - P(x, i)\|$$

This means that the current iteration is projected on one of the remotest sets of the family. Define the mapping $T: \mathcal{H} \to \mathcal{H}$ by $T(x) = P(x, i_x)$. It is clear that $x \in \cap M_i$ if and only if $T(x) = x$, hence if and only if $x$ is a fixed point of $T$, that is $\cap M_i = Fix(T)$. For any $x \in \mathcal{H}$ and $x^* \in Fix(T)$, the following Kolmogorov condition $\langle x - P(x, i_x), P(x, i_x) - x^* \rangle \geq 0$ is satisfied and it is routine to see that $T$ is demicontractive and we get

$$\|T_t(x) - x^*\|^2 < \|x - x^*\|^2 - t(2-t)\|x - T(x)\|^2 \qquad (5)$$

Therefore, $T_t$ is quasi-nonexpansive. According to Theorem (1), the sequence $\{x_k\}$ given by the generation function $T_t$ converges strongly to an element of $Fix(T)$ if and only if $d(x_k, Fix(T)) \to 0$ as $k \to \infty$

On the other hand, it is easy to see that $d(x_k, M_i) \to 0$ for each $i$. Indeed, from quasi-nonexpansivity of $T_t$ it follows that the sequence $\{\|x_k - y\|\}$ is monotone decreasing and bounded, therefore $\|x_k - y\| \to \delta_y$ as $k \to \infty$, for each $y \in \cap M_i$. From (5) we obtain that

$$\|x_k - T(x_k)\|^2 \leq \frac{1}{t(2-t)}(\|x_k - y\|^2 - \|x_{k+1} - y\|^2)$$

and hence $\|x_k - T(x_k)\| \to 0$ as $k \to \infty$. But $\|x - P(x,i)\| \leq \|x - T(x)\|$ for each $i$. Therefore $d(x_k, M_i) = \|x_k - P(x_k, i)\| \to 0$ as $k \to \infty$ and the essential point in the convergence of the projection method is the following property of the family $M_i$: *For any sequence $\{x_k\}$ such that $d(x_k, M_i) \to 0$ for each i, one has $d(x_k, \cap M_i) \to 0$.* Note that this property does not holds for any family (see the example below).

The above question was formulate by Gurin, Poliac and Raic [9] in 1967 in connection with strong convergence of the projection method. They proved that if $M_{\bar{\alpha}} \cap (Int \bigcap_{\alpha \in A} M_\alpha) \neq \emptyset$, where $M_{\bar{\alpha}}$ is a certain set of the family, then the family has the above property for any *bounded* sequence. In the sequel, we say that such a family has the GPR (Gurin, Poiac and Raic) property.

Bauschke and Borwein [2] introduced the notion of *regularity* for a finite family (N-tuple) of closed convex sets $M_1,...,M_N$ with nonempty intersection $M$, by the condition

$$\forall \varepsilon > 0, \quad \exists \delta > 0, \quad \forall x \in \mathcal{H}$$
$$max\{d(x, M_i), i = 1, \ldots, N\} \leq \delta$$
$$\Rightarrow d(x, M) \leq \varepsilon$$

If this holds only on bounded sets, then they speak of a boundedly regular family.

## 2 The regularity properties

In the paper [18] is pointed out that the GPR property of a family is in connection with the following two characteristics: the emptiness and the boundedness of the intersection of family. Therefore, the following four cases were

(1) *Int* $\cap$ $M_i \neq \emptyset$ and $\cap M_i$ is bounded;

(2) *Int* $\cap$ $M_i = \emptyset$ and $\cap M_i$ is unbounded,

(3) *Int* $\cap$ $M_i \neq \emptyset$ and $\cap M_i$ is unbounded;

(4) *Int* $\cap$ $M_i = \emptyset$ and $\cap M_i$ is bounded.

In [18] the first 3 cases were analyzed; for the case 1 it is proved that the GPR property holds, whereas for the cases 2 and 3 the property does not

231

holds (specific examples are given). The case 4 was left unsolved, the only remark that in this case the authors are looking for a suitable example is made. In the sequel we will prove that the GPR property holds for a particular family of sets with empty and bounded intersection (i.e. the case 4); so we give a partial answer for this case.

First of all we reproduce here the proof for the case 1 and the examples for the cases 2 and 3 from [17, 18].

**Case (1).** For the case of a bounded sequence, a similar result was given in [9].

**Lemma 1.** *Let $M_i \subset \mathcal{H}$ (i=1,...,m) be a family of convex sets such that Int $\cap M_i$ is nonempty and bounded and let $\{x_k\}$ be a sequence of $\mathcal{H}$ such that $d(x_k, M_i) \to 0$ as $k \to \infty$ for each i. Then $d(x_k, \cap M_i) \to 0$ as $k \to \infty$.*

*Proof.* We assume that $o \in Int \cap M_i$. Then there exists a closed ball $D(o,r)= \{x \in \mathcal{H} : \|x\| \le r\} \subset \cap M_i$. Let $\varepsilon$ be a given real number, $0 < \varepsilon < 1$, and let $C$ be a constant such that $\|x\| \le C-1$ for all $x \in \cap M_i$ which is possible, because $\cap M_i$ is bounded.

Since, $d(x_k, M_i) \to 0$ as $k \to \infty$ for each index $i$, there exists a sequence $\{y_k^{(i)}\}_{k \in N} \subset M_i$ such that $\|y_k^{(i)} - x_k\| \to 0$ as $k \to \infty$. Let

$$z_k = (1 - \tfrac{C}{\varepsilon})(y_k^{(i)} - x_k), \quad k = 0,1,\ldots \qquad (6)$$

There exists a number $k_i(\varepsilon)$ such that if $k \ge k_i(\varepsilon)$ then $\|y_k^{(i)} - x_k\| \le \tfrac{r}{|1-\tfrac{C}{\varepsilon}|}$ and so $\|z_k\| < r$, that is $z_k \in \cap M_i$.

On the other hand, from (6) we obtain

$$(1 - \tfrac{\varepsilon}{C})x_k = \tfrac{\varepsilon}{C}z_k + (1 - \tfrac{\varepsilon}{C})y_k^{(i)}$$

and for $k \ge k_i(\varepsilon)$ we have $(1-\tfrac{\varepsilon}{C})x_k \in M_i$, because $y_k^{(i)}, z_k \in M_i$ and $M_i$ are convex.

Now, let $k_0(\varepsilon)=max_i k_i(\varepsilon)$. Then, for $k \ge k_0(\varepsilon)$ it follows that $(1-\tfrac{\varepsilon}{C})x_k \in \cap M_i$ and

$$d(x_k, \cap M_i) \le \|x_k - (1-\tfrac{\varepsilon}{C})x_k\| = \tfrac{\varepsilon}{C-\varepsilon}\|(1-\tfrac{\varepsilon}{C})x_k\| < \varepsilon$$

which end the proof. ∎

**Case (2).** Example. Suppose that $\mathcal{H}$ is the real three-dimensional space, that the set $A_1$ is a cone and the set $A_2$ is a tangent plane $(ABCD)$. The situation is depicted in Figure 1

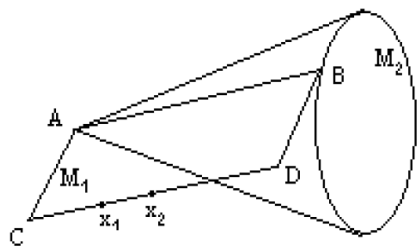

**Fig. 1.**

The plane $(ABCD)$ is tangent to the cone along the generatrix $(AB)$ and hence $A_1 \cap A_2 = (AB)$. Now, let us consider a sequence $\{x_k\}$ in the plane $(ABCD)$ such that $d(x_k, (AB)) = \delta = const.$ and $\|x_k\| \to \infty$ as $k \to \infty$. It is clear that $d(x_k, A_2) \to 0$ as $k \to \infty$ and $d(x_k, A_1) = 0$ for all $k$; but $d(x_k, A_1 \cap A_2) = \delta > 0$. Therefore, the conclusion of Lemma 1 is not true.

**Case (3).** The example is similar to that of the case (2). Suppose that $\mathcal{H}$ is the real three-dimensional space, that the family consists only of two sets, that $A_1$ is a cone and that $A_2$ is the half space defined by a secant plane parallel with a generatrix AB of the cone. The situation is depicted in Figure 2.

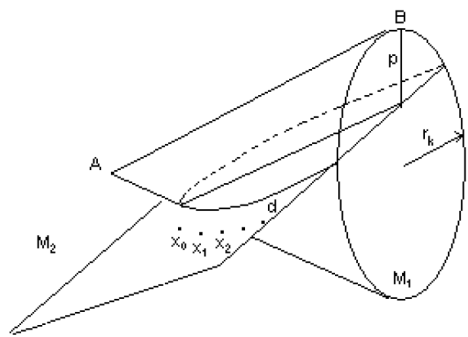

**Fig. 2.**

Obviously, $A_1 \cap A_2 \neq \emptyset$. Now, let us consider a sequence $\{x_k\}$ in the secant plane, with constant distances to the intersection parabola between cone and plane. We have $d(x_k, A_2)=0$. By an elementary computation, the distance between the terms of the sequence and $M_1$, is given by the formula

$$d(x_k, M_1) = \sqrt{r_k^2 + d^2 + 2d\sqrt{2pr_k - p^2}} - r_k,$$

where $p$, $d$, $r_k$ have the meaning from Fig. 1.

Therefore, $d(x_k, A_1) \to 0$, whereas $d(x_k, A_1 \cap A_2) = d > 0$.

**Case (4)**. We will prove that for a particular case the GPR property holds. This particular case is: The Hilbert space is the Euclidean 3 dimensional space, the family consists in two closed convex sets, $A$ and $B$, and the intersection $A \cap B$ is bounded and belongs to a plan $\mathcal{H}$ (thus the interior of the intersection is empty).

*Some preliminaries*. Let $A$ be a closed convex set in $\mathbb{R}^3$ and let $\mathcal{E}(A,x)$ be the convex cone hull of $A$ with vertex $x$. Let $d \in A$ and let $r = x+t(d-x)$ be the ray with origin $x$ passing though $d$. Let also $t_d$ be a positive number defined by

$$\|x + t_d(d - x)\| = \min_t \|x + t(d - x)\|$$

**Definition 5**. The superior side of $\mathcal{E}(A,x)$ is

$Sup\ \mathcal{E}(A,x) = \{x+t(d-x): 0 \leq t \leq t_d, d \in A\}$

**Lemma 2.** Let $A$ be a closed convex set in $\mathbb{R}^3$, let $x$ be a point such that $x \notin A$ and let $y$ be a point on the border of $Sup\ \mathcal{E}(A,x)$. If $P_x$ and $P_y$ denote the projection of $x$ and $y$ onto $A$ respectively, then

$$\|y - P_y\| \leq \|x - P_x\|$$

*Proof*. Let $z$ be a point in $Ray(x,y) \cap A$. Because $y$ belongs to the segment $[x,z]$, it follows that $y = tx+(1-t)z$, $0 \leq t \leq 1$. Let $u$ on be $[x,z]$ given by $y = tP_x + (1-t)z$; as $x$ and $P_x$ are in $A$, it follows that $u \in A$. From Kolmogorow characterization of the projection it obtain $\langle x - P_x, P_x - z \rangle \geq 0$. On the other hand, by a direct computation, it results

$$\|y - u\|^2 - \|y - P_y\|^2 \leq 2\langle y - P_y, P_y - u \rangle + \|P_y - u\|^2.$$

and $\|y - P_y\| \leq \|y - u\|$. Therefore,

$$\|y - P - y\| \leq \|y - u\| = t\|x - P_x\| \leq \|x - P_x\| \quad \blacksquare$$

Consider now the following enlargement of both sets $A$ and $B$: Let $x$ be a point does not belonging to $A$, and let $P_x$ be the projection of $x$ onto $A \cap B \in \mathcal{H}$. Suppose that $P_x \in Int(A \cap B)$. The symmetric point $x_s$ of $x$ with respect to $P_x$ lay in the opposite side of $\mathcal{H}$, therefore $x_s \notin B$.

Let $\mathcal{E}(A,x)$ and $\mathcal{E}(B,x_s)$ be the convex cone hulls of $A$ and $B$ with the vertex $x$ and $x_s$ respectively; these cone hulls are our enlargements.

**Lemma 3.** *The intersection $\mathcal{E}(A,x) \cap \mathcal{E}(B,x_s)$ is bounded and has its interior nonempty.*

*Proof.* Let $P_x r$ be a ray in $\mathcal{H}$ with the origin $P_x$. Because $A \cap B$ is bounded, there is a point $y \in P_x r$ such that $y \notin A \cap B$; for instance, suppose that $y \notin A$. Let $n$ be the normal line onto $\mathcal{H}$ in $y$. If for any $z \in n$, one has $z \in A$, then $y \in Fr(A)$ and with $y \notin A$, it follows that $A$ is an open set, which contradicts the hypothesis. Therefore, there is an $z \in n$ such that $z \notin A$ and also the ray $xz$ does not belongs to $\mathcal{E}(A,x)$.

On the other hand, let $x_s r$ be the external ray of $\mathcal{E}(B,x_s)$ in the plane $xP_x r$; in the worst case, this ray is parallel with $P_x r$. The construction is depicted in Figure 3.

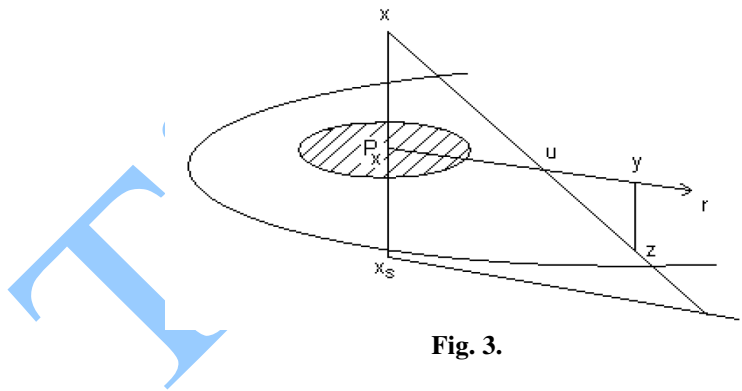

**Fig. 3.**

It is obvious that the intersection of $\mathcal{E}(A,x) \cap \mathcal{E}(B,x_s)$ in the plane $xP_x r$ is bounded. As the ray $P_x r$ was arbitrary, the same boundedness is valid for any ray with the origin $P_x$ and it results the boundedness of the intersection.

The fact that $Int\,\mathcal{E}(A,x) \cap \mathcal{E}(B,x_s)$ is nonempty, is obvious: a ball centered in $P_x$ sufficiently small belongs to $\mathcal{E}(A, x) \cap \mathcal{E}(B, x_s)$. ∎

We now are able to give a positive answer to the case 4.

**Lemma 4**. *Let A and B two closed convex sets in $\mathbb{R}^3$ so that $A \cap B$ is bounded and has its interior empty. Then the GPR property is valid.*

*Proof.* Obviously, $d(x_k, \mathcal{E}(A,x))$ and $d(x_k, \mathcal{E}(B,x_s))$ tend to zero; apply lemmas 1 and 3 to conclude that $d(x_k, \mathcal{E}(A,x) \cap \mathcal{E}(B,x_s)) \to 0$ So, given ε this distance becomes less then ε/2. We can take $x$ in lemma 3 such that $\|x - P_x\| \leq \varepsilon/2$. It follows that $d(x_k, A \cap B) \leq \varepsilon/2 + \varepsilon/2 = \varepsilon$. ∎